\documentclass[a4paper,10pt]{article}
\usepackage{color}   %May be necessary if you want to color links
\usepackage[bookmarks=true,pdfencoding=unicode]{hyperref}
\hypersetup{
    backref=true,
    citecolor=blue,
    colorlinks=true, %set true if you want colored links
    linktoc=all,     %set to all if you want both sections and subsections linked
    linkcolor=blue,  %choose some color if you want links to stand out
}

\usepackage{bookmark}
\usepackage[latin1]{inputenc}
\usepackage[T1]{fontenc}
\usepackage{bezier}
\usepackage{dsfont}
\usepackage{mathptmx}      % use Times fonts if available on your TeX system
\usepackage{amsfonts}
\usepackage{amsmath}
\usepackage{amssymb}
\usepackage{amsbsy}
\usepackage{amscd}
\usepackage{amsopn}
\usepackage{amssymb}
\usepackage{amstext}

\usepackage{amsxtra}
\usepackage{array}
\usepackage{bm}
\usepackage{graphicx,psfrag}
\usepackage{latexsym}
\usepackage{multirow}
\usepackage{subfig}
\usepackage{etoolbox}
\usepackage{algorithmicx}
\usepackage{algpseudocode}
\algtext*{EndIf}% Remove "end if" text
\usepackage{mdframed}
\usepackage{lipsum} % for creating dummy text

\newcommand{\wqed}{\hfill \ensuremath{\Box}}

\DeclareMathAlphabet{\mathcal}{OMS}{cmsy}{m}{n}

\newtoggle{LatexWithProofs}
\togglefalse{LatexWithProofs}

\newtoggle{LatexFull}
\togglefalse{LatexFull}

\newtheorem{pLemma}{Lemma}
\newtheorem{pProp}{Proposition}
\newtheorem{pArgument}{Argument}
\newtheorem{pExample}{Example}
\newtheorem{pCorol}{Corollary}
\newtheorem{pLongDef}{Definition}
\newtheorem{pTheorem}{Theorem}

\newcommand{\pRef}[1]{(\ref{#1})}

\newcommand{\pbFigureBH}[2]{
\begin{figure}[h]
\label{#1}
\hypertarget{#1}{}
\bookmark[
rellevel=1,
keeplevel,
dest=#1
]{Fig \ref{#1}: {#2}}
}

\newcommand{\peFigure}[1]{
\label{#1}
\end{figure}}

\newcommand{\peLongDef}[1]{$\blacktriangle$ \end{pLongDef}}
\newcommand{\peLongDefX}[1]{\end{pLongDef}}

\newcommand{\peArgument}[1]{ \hyperlink{#1}{$\blacktriangle$} \end{pArgument}}

\newcommand{\peExample}[1]{ \hyperlink{#1}{$\blacktriangle$} \end{pExample}}
\newcommand{\peExampleX}[1]{\end{pExample}}

\newcommand{\pePlainExample}[1]{$\blacktriangle$ \end{pExample}}
\newcommand{\pePlainExampleX}[1]{\end{pExample}}

\newcommand{\peCorol}[1]{\hyperlink{#1}{$\blacktriangle$} \end{pCorol}}
\newcommand{\peCorolX}[1]{\end{pCorol}}

\newcommand{\peLemma}[1]{ \hyperlink{#1}{$\blacktriangle$} \end{pLemma}}
\newcommand{\peLemmaX}[1]{\end{pLemma}}

\newcommand{\peProp}[1]{ \hyperlink{#1}{$\blacktriangle$} \end{pProp}}
\newcommand{\pePropX}[1]{\end{pProp}}

\newcommand{\peTheorem}[1]{ \hyperlink{#1}{$\blacktriangle$} \end{pTheorem}}

\iftoggle{LatexFull}
{
\newcommand{\pFullLink}[1]{\eqno \hyperlink{#1}{\blacktriangle}}
\newcommand{\peFullProp}[1]{ \hyperlink{#1}{$\blacktriangle$} \end{pProp}}
\newcommand{\peFullCorol}[1]{ \hyperlink{#1}{$\blacktriangle$} \end{pCorol}}
\newcommand{\peFullLemma}[1]{ \hyperlink{#1}{$\blacktriangle$} \end{pLemma}}
\newcommand{\peFullExample}[1]{ \hyperlink{#1}{$\blacktriangle$} \end{pExample}}
}
{
\newcommand{\pFullLink}[1]{\eqno \blacktriangle}
\newcommand{\peFullProp}[1]{$\blacktriangle$ \end{pProp}}
\newcommand{\peFullCorol}[1]{$\blacktriangle$ \end{pCorol}}
\newcommand{\peFullLemma}[1]{$\blacktriangle$ \end{pLemma}}
\newcommand{\peFullExample}[1]{$\blacktriangle$ \end{pExample}}
}

\iftoggle{LatexWithProofs}
{
\definecolor{darkpink}{rgb}{0.91, 0.33, 0.5}
\definecolor{darksalmon}{rgb}{0.91, 0.59, 0.48}
\definecolor{desertsand}{rgb}{0.93, 0.79, 0.69}
\definecolor{celadon}{rgb}{0.67, 0.88, 0.69}
\definecolor{darkcyan}{rgb}{0.0, 0.55, 0.55}

\newcommand{\pbClaim}[1]{{\color{red} ?}}
\newcommand{\peClaim}[1]{{\color{red} ?`}}
\newcommand{\pbArgument}[1]{\begin{pArgument} \label{#1}  {\color{darksalmon} #1}}

\newcommand{\pbExample}[1]{\begin{pExample} \label{#1}  {\color{darksalmon} #1}}

\newcommand{\pbExampleB}[1]{
\begin{pExample} \label{#1}
{\color{darksalmon} #1}
\hypertarget{{T#1}}{}
\bookmark[
rellevel=1,
keeplevel,
dest=T#1
]{Example \ref{#1}}
}

\newcommand{\pbExampleBT}[2]{
\begin{pExample}[#2] \label{#1}
{\color{darksalmon} #1}
\hypertarget{{T#1}}{}
\bookmark[
rellevel=1,
keeplevel,
dest=T#1
]{Example \ref{#1}: {#2}}
}

\newcommand{\pbCorol}[1]{\begin{pCorol} \label{#1}  {\color{darksalmon} #1}}

\newcommand{\pbCorolB}[1]{
\begin{pCorol} \label{#1}  {\color{darksalmon} #1}
\hypertarget{{T #1}}{}
\bookmark[
rellevel=1,
keeplevel,
dest=T#1
]{Corollary \ref{#1}}
}

\newcommand{\pbCorolBT}[2]{
\begin{pCorol}[#2]\label{#1}  {\color{darksalmon} #1}
\hypertarget{{T#1}}{}
\bookmark[
rellevel=1,
keeplevel,
dest=T#1
]{Cor. \ref{#1}: {#2}}
}

\newcommand{\pbLemma}[1]{\begin{pLemma} \label{#1}  {\color{darksalmon} #1}}

\newcommand{\pbLemmaB}[1]{
\begin{pLemma} \label{#1}  {\color{darksalmon} #1}
\hypertarget{{T#1}}{}
\bookmark[
rellevel=1,
keeplevel,
dest=T#1
]{Lemma \ref{#1}}
}

\newcommand{\pbLemmaBT}[2]{
\begin{pLemma}[#2] \label{#1}  {\color{darksalmon} #1}
\hypertarget{{T#1}}{}
\bookmark[
rellevel=1,
keeplevel,
dest=T#1
]{Lem. \ref{#1}: {#2}}
}

\newcommand{\pbProp}[1]{\begin{pProp} \label{#1}  {\color{darksalmon} #1}}

\newcommand{\pbPropB}[1]{
\begin{pProp} \label{#1}  {\color{darksalmon} #1}
\hypertarget{{T#1}}{}
\bookmark[
rellevel=1,
keeplevel,
dest=T#1
]{Prop. \ref{#1}}
}

\newcommand{\pbPropBT}[2]{
\begin{pProp}[#2] \label{#1}  {\color{darksalmon} #1}
\hypertarget{{T#1}}{}
\bookmark[
rellevel=1,
keeplevel,
dest=T#1
]{Prop. \ref{#1}: {#2}}
}

\newcommand{\pbTheorem}[1]{\begin{pTheorem} \label{#1}  {\color{darksalmon} #1}}

\newcommand{\pbTheoremB}[1]{
\begin{pTheorem} \label{#1}  {\color{darksalmon} #1}
\hypertarget{{T#1}}{}
\bookmark[
rellevel=1,
keeplevel,
dest=T#1
]{Theorem \ref{#1}}
}

\newcommand{\pbTheoremBT}[2]{
\begin{pTheorem}[#2] \label{#1}  {\color{darksalmon} #1}
\hypertarget{{T#1}}{}
\bookmark[
rellevel=1,
keeplevel,
dest=T#1
]{Theor. \ref{#1}: {#2}}
}

\newcommand{\pbLongDef}[1]{\begin{pLongDef}
\label{#1} \hypertarget{{#1}}{} {\color{darksalmon} #1}
\bookmark[
rellevel=1,
keeplevel,
dest=#1
]{Definition \ref{#1}}
}

\newcommand{\pbLongDefB}[2]{\begin{pLongDef}
\label{#1} \hypertarget{{#1}}{} {\color{darksalmon} #1}
\bookmark[
rellevel=1,
keeplevel,
dest=#1
]{\ref{#1}: #2}
}

\newcommand{\pbLongDefBT}[2]{\begin{pLongDef}[#2]
\label{#1} \hypertarget{{#1}}{} {\color{darksalmon} #1}
\bookmark[
rellevel=1,
keeplevel,
dest=#1
]{\ref{#1}: #2}
}

\newcommand{\pbChain}[1]{\[ }
\newcommand{\peChain}[1]{\ {\color{red} \checkmark  \wrm{#1}} \] }

\newcommand{\pePClaim}[1]{ \checkmark}

\newcommand{\pbTClaim}[1]{\begin{equation} \label{#1}  }
\newcommand{\peTClaim}[1]{ \ {\color{red} \checkmark  \wrm{#1}} \end{equation} }
\newcommand{\peEClaim}[1]{ \ {\color{red} \checkmark  \wrm{#1 \bullet} }\end{equation} }

\newcommand{\pbDef}[1] {\hypertarget{{#1}}{} \begin{equation} \label{#1} }
\newcommand{\peDef}[1]{ {\color{blue} \bigstar \wrm{#1}} \end{equation}}

\newcommand{\pbHypot}[1]{ \begin{equation} \label{#1}  }
\newcommand{\peHypot}[1]{  \ \ {{\color{blue} \bigstar \wrm{#1}}}  \end{equation} }

\newcommand{\pbProofB}[2]{\newpage {\color{darksalmon} Proof of {#1} {#2}.} \ref{#2} \hypertarget{{#2}}{}
\bookmark[
rellevel=1,
keeplevel,
dest=#2
]
{{#1} \ref{#2}}
}

\newcommand{\peProof}[2]{\wqed{} {\color{darksalmon} End of proof of #1 #2} \ref{#2}}

\newcommand{\peVerify}[1]{\wqed{} {\color{darksalmon} End of verification of Example #1} \ref{#1}}

}
{
\newcommand{\pbChain}[1]{\[}
\newcommand{\peChain}[1]{\]}

\newcommand{\pbClaim}[1]{}
\newcommand{\peClaim}[1]{}

\newcommand{\peTClaim}[1]{ \end{equation}}
\newcommand{\peEClaim}[1]{ \end{equation}}
\newcommand{\pbTClaim}[1]{\begin{equation} \label{#1}}

\newcommand{\peDef}[1]{  \end{equation}}

\newcommand{\lpPlainClaim}[1]{}
\newcommand{\lpEndPlainClaim}[1]{}
\newcommand{\pbDef}[1]{ \hypertarget{{#1}}{} \begin{equation} \label{#1}}

\newcommand{\pbHypot}[1]{\begin{equation} \label{#1}  }
\newcommand{\peHypot}[1]{ \end{equation}}

\newcommand{\pbCorol}[1]{\begin{pCorol} \label{#1}}

\newcommand{\pbCorolB}[1]{\begin{pCorol} \label{#1}
\hypertarget{{T#1}}{}
\bookmark[
rellevel=1,
keeplevel,
dest=T#1
]
{Corollary \ref{#1}}
}

\newcommand{\pbCorolBT}[2]{
\begin{pCorol}[#2] \label{#1}
\hypertarget{{T#1}}{}
\bookmark[
rellevel=1,
keeplevel,
dest=T#1
]
{Cor. \ref{#1}: {#2}}
}

\newcommand{\pbArgument}[1]{\begin{pArgument} \label{#1}}
\newcommand{\pbExample}[1]{\begin{pExample} \label{#1}}

\newcommand{\pbExampleB}[1]{\begin{pExample} \label{#1}
\hypertarget{{T#1}}{}
\bookmark[
rellevel=1,
keeplevel,
dest=T#1
]{Example \ref{#1}}
}

\newcommand{\pbExampleBT}[2]{
\begin{pExample}[#2] \label{#1}
\hypertarget{{T#1}}{}
\bookmark[
rellevel=1,
keeplevel,
dest=T#1
]{Example \ref{#1}: {#2}}
}

\newcommand{\pbLemma}[1]{\begin{pLemma} \label{#1}}
\newcommand{\pbLemmaBT}[2]{\begin{pLemma}[#2] \label{#1}
\hypertarget{{T#1}}{}
\bookmark[
rellevel=1,
keeplevel,
dest=T#1
]
{Lem. \ref{#1}: {#2}}
}

\newcommand{\pbLemmaB}[1]{\begin{pLemma} \label{#1}
\hypertarget{{T#1}}{}
\bookmark[
rellevel=1,
keeplevel,
dest=T#1
]
{Lemma \ref{#1}}
}

\newcommand{\pbProp}[1]{\begin{pProp} \label{#1}}
\newcommand{\pbPropBT}[2]{\begin{pProp}[#2] \label{#1}
\hypertarget{{T#1}}{}
\bookmark[
rellevel=1,
keeplevel,
dest=T#1
]
{Prop. \ref{#1}: {#2}}
}

\newcommand{\pbPropB}[1]{\begin{pProp} \label{#1}
\hypertarget{{T#1}}{}
\bookmark[
rellevel=1,
keeplevel,
dest=T#1
]
{Prop. \ref{#1}}
}

\newcommand{\pbTheorem}[1]{\begin{pTheorem} \label{#1}}

\newcommand{\pbTheoremBT}[2]{\begin{pTheorem}[#2] \label{#1}
\hypertarget{{T#1}}{}
\bookmark[
rellevel=1,
keeplevel,
dest=T#1
]
{Thm. \ref{#1}: {#2}}
}

\newcommand{\pbLongDef}[1]{\begin{pLongDef} \label{#1}}

\newcommand{\pbLongDefB}[2]{\begin{pLongDef}
\label{#1} \hypertarget{{#1}}{}
\bookmark[
rellevel=1,
keeplevel,
dest=#1
]{\ref{#1}: #2}
}

\newcommand{\pbLongDefBT}[2]{\begin{pLongDef}[#2]
\label{#1} \hypertarget{{#1}}{}
\bookmark[
rellevel=1,
keeplevel,
dest=#1
]{\ref{#1}: #2}
}

\newcommand{\pbProofB}[2]{ {\bf Proof of {#1} \ref{#2}.} \hypertarget{{#2}}{}
\bookmark[
rellevel=1,
keeplevel,
dest=#2
]
{{#1} \ref{#2}}
}

\newcommand{\peProof}[2]{\wqed{}}

\newcommand{\peVerify}[1]{\wqed{}}
}

\newcommand{\wabs}[1]{\left|#1\right|}

\newcommand{\wc}[1]{\mathds{C}}
\newcommand{\wcal}[1]{\mathcal{#1}}

\newcommand{\whessf}[2]{{\nabla^2 \! #1 \left( #2 \right)}}

\newcommand{\wfc}[2]{{#1}\!\left(#2\right)}

\newcommand{\wgradf}[2]{{\nabla \! #1}\left(#2\right)}

\newcommand{\wi}[1]{\wrm{i}}

\algdef{SE}[DOWHILE]{Do}{doWhile}{\algorithmicdo}[1]{\algorithmicwhile\ #1}%

\newdimen\CdotAxis
\newcommand*{\CdotAux}[3]{%
  {%
    \settoheight\CdotAxis{$#2\vcenter{}$}%
    \sbox0{%
      \raisebox\CdotAxis{%
        \scalebox{#1}{%
          \raisebox{-\CdotAxis}{%
            $\mathsurround=0pt #2#3$%
          }%
        }%
      }%
    }%
    \dp0=0pt %
    \sbox2{$#2\bullet$}%
    \ifdim\ht2<\ht0 %
      \ht0=\ht2 %
    \fi
    \sbox2{$\mathsurround=0pt #2#3$}%
    \hbox to \wd2{\hss\usebox{0}\hss}%
  }%
}
\newcommand{\wlr}[1]{\left( #1 \right)}

\newcommand{\wn}{\mathds N}

\newcommand{\wnorm}[1]{\left\| #1 \right\|}

\newcommand{\wrone}{\mathds R}
\newcommand{\wrn}[1]{{\mathds R}^{#1}}
\newcommand{\wrm}[1]{\mathrm{#1}}

\newcommand{\wset}[1]{{\left\{ #1 \right\}}}

\newcommand{\wtr}{\mathrm{T}}

\newcommand{\wvec}[1]{\mathbf{#1}}

\newcommand{\wfpes}[1]{\wcal{E}_{\wcal{A},e}}

\newcolumntype{M}[1]{>{\centering\arraybackslash}m{#1}}
\newcolumntype{N}{@{}m{0pt}@{}}
\begin{document}
\title{A simple canonical form for nonlinear programming problems and its use}

\author{Walter F. Mascarenhas\thanks{
Instituto de Matem\'{a}tica e Estat\'{i}stica, Universidade de S\~{a}o Paulo,
Cidade Universit\'{a}ria, Rua do Mat\~{a}o 1010, S\~{a}o Paulo SP, Brazil. CEP 05508-090.
Tel.: +55-11-3091 5411, Fax: +55-11-3091 6134, walter.mascarenhas@gmail.com.}}
\markboth{W.F. Mascarenhas}{Stability versus order of approximation in barycentric interpolation}
\maketitle

\begin{abstract}
We argue that reducing nonlinear programming problems to
a simple canonical form is an effective way to analyze
them, specially when the problem is degenerate and the
usual linear independence hypothesis does not hold.
To illustrate this fact we solve an open
problem about constraint qualifications using
 this simple canonical form.
\end{abstract}

\section{Introduction}
In this article we look at the classical nonlinear programming problem
\pbDef{nlp}
\begin{array}{llcl}
\wrm{minimize}      &  \wfc{f}{\wvec{x}} &      & \\
\wrm{subject \ to } &  \wfc{h}{\wvec{x}} & =    & 0 \in \wrn{m}, \\
                    &  \wfc{g}{\wvec{x}} & \leq & 0 \in \wrn{r + p}
\end{array}
\peDef{nlp}
in situations in which the derivatives $\wfc{Dh}{0}$
and $\wfc{Dg}{0}$ of the constraints are defective, in the
sense that the rank of combined matrix
\pbDef{jhg}
\wfc{Dhg}{\wvec{x}}
:=
\left(
\begin{array}{c}
\wfc{Dh}{\wvec{x}} \\
\wfc{Dg}{\wvec{x}}
\end{array}
\right)
\peDef{jhg}
is $m + r \leq m + r + p$ (by ``derivative'' $\wfc{Dh}{\wvec{x}}$
here we mean the linear transformation represented
by the Jacobian matrix of the function $h$.)
To simplify the notation,
we look at problem \pRef{nlp} for $\wvec{x}$ in an open
set $\wcal{A}$ containing $0 \in \wrn{n}$ and
assume that all constraints are active at $0$,
that is $\wfc{h}{0} = 0$ and $\wfc{g}{0} = 0$.

The analysis of problem
\pRef{nlp} is usually based on a {\it constraint qualification}.
For instance, the well known Mangasarian-Fromovitz constraint qualification
assumes that $\wfc{Dh}{0}$ has rank $m$ and require the
existence of a ``strictly decreasing direction $\wvec{d}$'' for
the inequality constraints, which is also compatible with the
equality constraints:

\pbLongDefBT{defMF}{The Mangasarian-Fromovitz Constraint Qualification}
For functions $h$ and $g$ such that
$\wfc{h}{0} = 0$ and $\wfc{g}{0} = 0$, we say that
$0 \in \wrn{n}$  satisfies the Mangasarian-Fromovitz constraint qualification for
problem
\pRef{nlp} if $\wfc{rank}{\wfc{Dh}{0}} = m$
and there exists $\wvec{d} \in \wrn{n}$ such that
$\wfc{Dh}{0} \wvec{d} = 0$ and
$\wfc{Dg}{0} \wvec{d} < 0$.
\peLongDef{defMF}

Unfortunately, even under conditions like Mangasarian-Fromovitz,
the analysis of problem \pRef{nlp} can be tricky, and it is
not uncommon to find articles in which incorrect results
or proofs are presented, as pointed out in \cite{Gabriel}.
In fact, we are all human and our capability to deal with the details involved
in the analysis of problem \pRef{nlp} is limited,
and we need tools to handle them.
In the present article we
argue that by reducing problem \pRef{nlp} to
a simple canonical form we have a better chance of understanding
degenerated cases, and we illustrate this point
by proving the following conjecture by Andreani, Martinez and Schuverdt
\cite{Nino} regarding second order constraint
qualifications.

\pbTheoremBT{thmNino}{Andreani's Conjecture}
Suppose the functions $f$, $h$ and $g$ in problem \pRef{nlp} are
of class $C^2$ in a neighborhood $\wcal{A}$ of $0 \in \wrn{n}$
and $\wfc{h}{0} = 0$ and $\wfc{g}{0} = 0$.
If the Mangasarian-Fromovitz constraint qualification is satisfied
and
\[
\wfc{\wrm{rank}}{\wfc{Dhg}{\wvec{x}}} \leq
\wfc{\wrm{rank}}{\wfc{Dhg}{0}} + 1
\]
for $\wvec{x} \in \wcal{A}$ then there exist
$\lambda \in \wrn{m}$
and $\mu \in \wrn{p}$, with $\mu \geq 0$,
such that
\pbDef{eqWeak}
\wfc{S}{\wvec{d}} :=  \wvec{d}^{\wtr{}}
\wlr{
\whessf{f}{0} + \sum_{j = 1}^m \lambda_j \whessf{h_j}{0}
+ \sum_{\ell = 1}^{r + p} \mu_\ell \whessf{g_\ell}{0}
}
\wvec{d} \geq 0
\peDef{eqWeak}
for all $\wvec{d} \in \wrn{n}$ with $\wfc{Dh}{0} \wvec{d} = 0$ and
$\wfc{Dg}{0} \wvec{d} = 0$.
\peTheorem{thmNino}

The reduction of complicated problems to simpler ones is a standard
procedure in mathematics. Frequently, the Inverse Function Theorem and
its variations, like the Implicit Function Theorem, are used
to obtain changes of variables which reduce a nonlinear
problem to a linear one. The results can be striking, as in
Malgrange's Preparation Theorem, and Mather's Division Theorem
\cite{Golubitsky}, which generalize the Implicit Function Theorem
and have their roots in Weierstrass' work.
The Hartman-Grobman Theorem \cite{Perko,Strogatz} is another
remarkable result along these lines, and allows us
to understand the behavior of a nonlinear dynamical system
in terms of its linearization, under mild technical conditions.

The Stable Manifold Theorem  \cite{Perko} is yet another example.
In our research about nonlinear programming we have
used it to analyze the convergence
of significant algorithms, like the Affine Scaling
Method \cite{Affine}, Newton's Method
\cite{Newton}, the BFGS Method \cite{BFGS},
and even to general families of methods
\cite{Divergence,GaussNewton}. In these articles,
by looking at low
dimensional problems from the right perspective,
we were able to provide satisfactory answers to
relevant open problems.

Here we present yet another instance in
which the proper simplification
leads to a better understanding of
nonlinear programming problems.
In Section \ref{secChange} we emphasize
that usual optimality conditions in
nonlinear programming are invariant under
changes of variables. In particular the
Lagrange multipliers are the same in
different coordinate systems for the
independent variables, and the Mangasiran-Fromovitz
constraint qualification is invariant under
changes of these variables.
Section \ref{secChange} is 
obvious but relevant. In fact, we believe that
changes of variables do not receive the attention
they deserve in the mathematical programming
literature. For instance,
most textbooks do not mention the invariance  above
explicitly. Authors usually have these
changes of variables on the back of their
mind, and build good examples to illustrate their points based
upon them.
However, readers may not notice that, with the proper changes
of variables, more general situations can be
reduced to these good examples, and in many
cases such examples are more enlightening than proofs.

In Section \ref{secRank} we present a convenient
change of the independent variables for nonlinear programming
problems in which the rank of the derivative
of the inequality constraints drops at the point in
which we are interested. Of course, we
cannot rely on the usual linear independence
hypothesis in this case, but we show that
we can still simplify problem considerably
by using the proper choice of variables.
Finally, in Section \ref{secNino} we prove
Andreani's conjecture, using the results from
the previous sections and two linear
algebraic lemmas proved in the appendix.

\section{Changes of variables}
\label{secChange}
This section calls the reader's attention to
the fundamental fact that the usual optimality conditions
in mathematical programming are invariant
under changes of variables. In other
words,  in theory we can pick any coordinate
system we please, as long as the coordinates
are ``consistent'' and reflect correctly
the problem we want to understand.

We emphasize that, for theoretical
purposes, changes of coordinates do not
need to be explicit: it suffices to know that
they exist and have the nice properties required
to deal with the problem at hand.
Of course, things are different in practice,
because we usually cannot afford, or know
how, to compute such changes of coordinates.

The simplest changes of variables are the
linear ones: we replace the coordinates
$\wvec{x}$ by $\wvec{A} \wvec{y}$, where
$\wvec{A}$ is a nonsingular square matrix.
As a result, we can replace a function $f$
by a more appropriate, or simpler, function
$\hat{f}$, as in
\pbDef{fullLin}
\wfc{\hat{f}}{y} = \wfc{f}{\wvec{A} \wvec{y}}.
\peDef{fullLin}
The chain rule yields
%
%\pbDef{djh}
%\partial_j \wfc{\hat{f}}{\wvec{y}} = \sum_{k = 1}^m \partial_k \wfc{f}{\wvec{A} \wvec{y}} \, a_{k j} =
%\sum_{k = 1}^m  a_{j k}^{\wtr{}} \partial_k \wfc{f}{\wvec{A} \wvec{y}} =
%\wlr{\wvec{A}^{\wtr{}} \wgradf{f}{\wvec{A} \wvec{y}}}_j
%\peDef{djh}
%

%
\pbDef{linGrad}
\wgradf{\hat{f}}{\wvec{y}} = \wvec{A}^{\wtr{}} \wgradf{f}{\wvec{A} \wvec{y}}.
\peDef{linGrad}
and applying the formula above to the function
$\wfc{\hat{f}_k}{\wvec{y}} := \partial_k \wfc{f}{\wvec{A} \wvec{y}}$
and taking one more derivative
\if{false}
we obtain
\[
\partial_{i,j} \wfc{\hat{f}}{\wvec{y}} = \sum_{k = 1}^m \sum_{\ell = 0}^m \partial_{\ell,k} \wfc{f}{\wvec{A} \wvec{y}} \, a_{\ell,i} a_{k j} =
\sum_{\ell = 1}^m a_{\ell,i} \sum_{k = 1}^m \partial_{\ell,k} \wfc{f}{\wvec{A} \wvec{y}} \,  a_{k j} =
\]
\[
\sum_{\ell = 1}^m a_{\ell,i} \sum_{k = 1}^m \partial_{\ell,k} \wfc{f}{\wvec{A} \wvec{y}} \,  \wlr{\wvec{A} \wvec{e}_j}_k =
\sum_{\ell = 1}^m a_{\ell,i} \wlr{\whessf{f}{\wvec{A} \wvec{y}} \, \wvec{A} \wvec{e}_j}_\ell =
\]
\[
\sum_{\ell = 1}^m a_{i,\ell}^{\wtr{}} \wlr{\whessf{f}{\wvec{A} \wvec{y}} \, \wvec{A} \wvec{e}_j}_\ell =
\sum_{\ell = 1}^m \wlr{\wvec{e}_i^\wtr{} \wvec{A}}_\ell \wlr{\whessf{f}{\wvec{A} \wvec{y}} \, \wvec{A} \wvec{e}_j}_\ell
= \wvec{e}_i^{\wtr{}} \wvec{A}^{\wtr{}} \whessf{f}{\wvec{A} \wvec{y}} \, \wvec{A} \wvec{e_j},
\]
and this implies that
\fi
and doing the algebra we obtain that
\pbDef{linHess}
\whessf{\hat{f}}{\wvec{y}} = \wvec{A}^{\wtr{}}\,  \whessf{f}{\wvec{A} \wvec{y}} \, \wvec{A}.
\peDef{linHess}
Equations \pRef{linGrad} and \pRef{linHess}
are quite useful, but they are not enough to explore
the full power of changes of variables. For that we
we need to replace $\wvec{x}$ by a nonlinear
function $\wfc{q}{\wvec{y}}$ of a more convenient
variable $\wvec{y}$, as in the nonlinear
version of Equation \pRef{fullLin}:
\pbDef{fullNL}
\wfc{\hat{f}}{\wvec{y}} := \wfc{f}{\wfc{q}{\wvec{y}}}.
\peDef{fullNL}
In Equation \pRef{fullNL},
the function $q$ is a {\it local diffeomorphism}, that is
a differentiable function defined in a neighborhood
of the point with which we are concerned, and such
that its inverse $q^{-1}$ (in the sense that
$\wfc{q^{-1}}{\wfc{q}{\wvec{x}}} = \wvec{x}$) exists and
is also differentiable.
In this nonlinear setting we have the following
version of Equation \pRef{linGrad}
\pbDef{nonGrad}
\wgradf{\hat{f}}{y} = \wfc{Dq}{\wvec{y}}^{\wtr{}} \wgradf{f}{\wfc{q}{\wvec{y}}}.
\peDef{nonGrad}
Equation \pRef{nonGrad} is almost the same as \pRef{linGrad}: we only need to
replace $\wvec{A}$ by $\wfc{Dq}{\wvec{y}}$ and $\wvec{A} \wvec{y}$ by
$\wfc{q}{\wvec{y}}$, and keep in mind that the matrix
$\wfc{Dq}{\wvec{y}}$ is square and non singular.
Things are a bit more complicated for the Hessian. In this case we need to
introduce an extra sum due to the curvature
in $q$, and the resulting formula is:
%
%Taking one more derivative of \pRef{djq} we obtain
%\[
%\partial_{i,j} \wfc{\hat{f}}{\wvec{y}} =
%\sum_{k = 1}^n \sum_{\ell = 0}^m \partial_{\ell, k} \wfc{f}{\wfc{q}{\wvec{y}}}
%\, \partial_i \wfc{q_\ell}{\wvec{y}} \, \partial_j \wfc{q_k}{\wvec{y}}
%+
%\sum_{k = 1}^n \partial_k \wfc{f}{\wfc{q}{\wvec{y}}} \, \partial_{i,j} \wfc{q_{k}}{\wvec{y}}.
%\]
%As before, it follows that
%%%
\pbDef{nonHess}
\whessf{\hat{f}}{\wvec{y}} = \wfc{Dq}{\wvec{y}}^{\wtr{}} \whessf{f}{\wfc{q}{\wvec{y}}} \wfc{Dq}{\wvec{y}}
+ \sum_{k=1}^n \partial_k \wfc{f}{\wfc{q}{\wvec{y}}} \whessf{q_k}{\wvec{y}}.
\peDef{nonHess}
%%%

With Equations \pRef{nonGrad} and \pRef{nonHess}
we can find how the optimality conditions behave under
nonlinear changes of coordinates $x = \wfc{q}{\wvec{y}}$.
To see  why this is true,
we consider the classical nonlinear
programming problem \pRef{nlp}.
By making the change of variables $\wvec{x} = \wfc{q}{\wvec{y}}$ we do
not affect the satisfiability of the equalities and inequalities
in \pRef{nlp}, that is
\[
\wfc{\hat{h}_j}{\wvec{y}} := \wfc{h_j}{\wfc{q}{\wvec{y}}}
\]
will be equal to zero as long as
$\wfc{h_j}{\wvec{x}}$ is equal to zero. This
is obvious, but the analogous obvious property
would not hold if we were to take combinations of the dependent
instead of the independent variables in problem \pRef{nlp},
as when replace the equations
\pbDef{depend}
\wfc{g_1}{\wvec{x}} \leq 0 \hspace{1cm} \wrm{and} \hspace{1cm}  \wfc{g_2}{\wvec{x}} \leq 0
\peDef{depend}
by
\[
\wfc{g_1}{\wvec{x}} + \wfc{g_2}{\wvec{x}} \leq 0 \hspace{1cm} \wrm{and}
\hspace{1cm}  \wfc{g_1}{\wvec{x}} + 2 \wfc{g_2}{\wvec{x}} \leq 0.
\]
Due to this obvious fact we must think carefully
before using techniques like the SVD to analyze the nonlinear
programming problem \pRef{nlp}. In particular,
when we say ``change of variables''
in this article we refer to changes in the dependent variables
$\wvec{x}$, but not to combinations of the equality constraints or
changes in the order of the inequality constraints.

Let us now analyze the first order optimality conditions for
problem \pRef{nlp}.
These conditions are written in terms of the Lagrange multiplies
$\lambda_j \in \wrone{}$ and $\mu_\ell \geq 0$:
\pbDef{firstOrder}
\wgradf{f}{\wvec{x}} +
\sum_{j = 1}^m \lambda_j \wgradf{q_j}{\wvec{x}} +
\sum_{\ell = 1}^{r + p} \mu_\ell \wgradf{g_\ell}{\wvec{x}} = 0,
\peDef{firstOrder}
Given a local diffeomorphism  $q$ with $\wfc{q}{\wvec{y}} = \wvec{x}$,
by defining
\pbDef{fgHat}
\wfc{\hat{f}}{\wvec{y}}                := \wfc{f}{\wfc{q}{\wvec{y}}},
\hspace{0.5cm} \wfc{\hat{h}}{\wvec{y}} := \wfc{h}{\wfc{q}{\wvec{y}}},
\hspace{0.5cm} \wrm{and} \hspace{0.5cm}
\wfc{\hat{g}}{\wvec{y}} := \wfc{g}{\wfc{q}{\wvec{y}}},
\peDef{fgHat}
Equation \pRef{nonGrad} yields
\pbDef{gradF}
\wgradf{\hat{f}}{\wvec{y}} = \wfc{Dq}{\wvec{y}}^{\wtr{}} \wgradf{f}{\wfc{q}{\wvec{y}}}
\peDef{gradF}
and
\pbDef{gradHat}
\wgradf{\hat{h}_j}{\wvec{y}} = \wfc{Dq}{\wvec{y}}^{\wtr{}} \wgradf{h_j}{\wfc{q}{\wvec{y}}}
\hspace{1cm} \wrm{and} \hspace{1cm}
\wgradf{\hat{g}_\ell}{\wvec{y}} = \wfc{Dq}{\wvec{y}}^{\wtr{}} \wgradf{g_\ell}{\wfc{q}{\wvec{y}}}.
\peDef{gradHat}
Since our $\wfc{Dq}{\wvec{y}}$ is always non singular,
it is easy to see that
the first order condition \pRef{firstOrder} holds if and only if
\pbDef{firstOrderY}
\wgradf{\hat{f}}{\wvec{y}}
+ \sum_{j = 1}^m \lambda_j \wgradf{\hat{h}_j}{\wvec{y}}
+ \sum_{\ell = 1}^{r + p} \mu_\ell \wgradf{\hat{g}_\ell}{\wvec{y}} = 0.
\peDef{firstOrderY}
In other words, the first order conditions are invariant with
respect to changes of coordinates, and we may study them
by using Equation \pRef{firstOrder},
or Equation \pRef{firstOrderY}, or both.
In particular, the Lagrange multipliers do not
change as we change coordinates as above.

For the second order conditions, we consider
directions $\wvec{d}$ such that
\pbDef{ortho}
\wfc{Dh}{0} \wvec{d}  = 0
\hspace{1cm} \wrm{and} \hspace{1cm}
\wfc{Dg}{0} \wvec{d} = 0
\peDef{ortho}
Equation \pRef{gradHat} shows that
Equation \pRef{ortho} is equivalent to
\pbDef{orthoY}
\wfc{D\hat{h}}{0} \hat{\wvec{d}} = 0
\hspace{1cm} \wrm{and} \hspace{1cm}
\wfc{D\hat{g}}{0} \hat{\wvec{d}} = 0,
\peDef{orthoY}
for
\[
\hat{\wvec{d}} := \wfc{Dq}{\wvec{y}}^{-1} \wvec{d},
\hspace{0.5cm} \wrm{or} \hspace{0.5cm}
\wvec{d} = \wfc{Dq}{\wvec{y}} \hat{\wvec{d}},
\]
and the orthogonality condition \pRef{ortho} is invariant
under changes of coordinates, that is, Equations
\pRef{ortho} for $\wvec{d}$ and Equation
\pRef{orthoY} for $\hat{\wvec{d}}$ are equivalent.

Using the equations above, we can write the second order term
\pbDef{termS}
\wfc{S}{\wvec{d}} := \wfc{S}{\wvec{d},f,h,g,\lambda, \mu} :=
\wvec{d}^{\wtr{}} \wlr{
\whessf{f}{\wvec{x}}  +
\sum_{j=1}^m \lambda_j \whessf{h_j}{\wvec{x}} \wvec{d} +
\sum_{\ell=1}^{r + p} \mu_\ell \whessf{g_\ell}{\wvec{x}}} \wvec{d},
\peDef{termS}
as
\[
\wfc{S}{\wvec{d}} =
\hat{\wvec{d}}^{\wtr{}} \left(
\wfc{Dq}{\wvec{y}}^{\wtr{}} \whessf{f}{\wfc{h}{\wvec{y}}} \wfc{Dq}{\wvec{y}} \ + \
\right.
\]
\[
\left.
\sum_{j=1}^m \lambda_j \wfc{Dq}{\wvec{y}}^{\wtr{}} \whessf{h_j}{\wfc{q}{\wvec{y}}} \wfc{Dq}{\wvec{y}} +
\sum_{\ell=1}^{r + p} \mu_\ell \wfc{Dq}{\wvec{y}}^{\wtr{}} \whessf{g_\ell}{\wfc{q}{\wvec{y}}} \wfc{Dq}{\wvec{y}}\right) \hat{\wvec{d}}.
\]
Equation \pRef{nonHess} yields
\[
\wfc{S}{\wvec{d},f,h,g,\lambda, \mu} =
\wfc{S}{\hat{\wvec{d}},\hat{f},\hat{h}, \hat{g},\lambda, \mu} - \Delta
\]
for
\[
\Delta
:=
\hat{\wvec{d}}^{\wtr{}} \left(
\sum_{k = 1}^n \partial_k \wfc{f}{\wfc{h}{\wvec{y}}} \whessf{q_k}{\wvec{y}} \ + \
\right.
\]
\[
\left.
\sum_{j = 1}^m \lambda_{j} \sum_{k=1}^n \partial_k \wfc{h_j}{\wfc{q}{\wvec{y}}} \whessf{q_k}{\wvec{y}} +
\sum_{\ell = 1}^{r + p} \mu_{\ell} \sum_{k=1}^n \partial_k \wfc{g_\ell}{\wfc{q}{\wvec{y}}} \whessf{q_k}{\wvec{y}}\right) \hat{\wvec{d}}
\]
\[
=  \hat{\wvec{d}}^{\wtr{}} \wlr{
\sum_{k=1}^n \wlr{\partial_k \wfc{f}{\wfc{q}{\wvec{y}}} +
\sum_{j = 1}^m \lambda_{j} \partial_k \wfc{h_j}{\wfc{q}{\wvec{y}}} +
\sum_{\ell = 1}^{r + p} \mu_{\ell} \partial_k \wfc{g_\ell}{\wfc{q}{\wvec{y}}} }
 \whessf{q_k}{\wvec{y}}} \hat{\wvec{d}}.
\]
When the first order conditions \pRef{firstOrder} hold we have that
\[
\partial_k \wfc{f}{\wfc{h}{\wvec{y}}}
+ \sum_{j = 1}^m \lambda_j \wgradf{h_j}{\wfc{q}{\wvec{y}}}
+ \sum_{\ell = 1}^{r + p} \mu_\ell \wgradf{g_\ell}{\wfc{q}{\wvec{y}}} = 0
\]
and $\Delta = 0$. Therefore, when the first order conditions hold, the
second order term $\wfc{S}{\wvec{d}}$ in Equation \pRef{termS}
is invariant with respect to changes of variables,
and it can be evaluated using the expression
\pbDef{termSY}
\wfc{S}{\wvec{d}} =
\wfc{\hat{S}}{\hat{\wvec{d}}} := \hat{\wvec{d}}^{\wtr{}} \wlr{
\whessf{\hat{f}}{\wvec{y}}  +
\sum_{j=1}^m \lambda_j \whessf{\hat{h}_j}{\wvec{y}} +
\sum_{\ell=1}^{r + p} \mu_\ell \whessf{\hat{g}_\ell}{\wvec{y}}} \hat{\wvec{d}}.
\peDef{termSY}

An analogous argument starting from the conditions
\[
\wfc{Dh}{0} \wvec{d}  = 0
\hspace{1cm} \wrm{and} \hspace{1cm}
\wfc{Dg}{0} \wvec{d} < 0
\]
instead of Equation \pRef{ortho} shows that the Mangazarian-Fromovitz
constraint qualification is invariant under changes of coordinates,
and similar arguments apply to many other constraint qualifications.

In summary, when trying to answer many theoretical questions regarding
the first and second order optimality conditions, Lagrange
multipliers and constraint qualifications for the nonlinear programming
problem \pRef{nlp}, we can analyze them in other coordinate systems, and
reach correct conclusions by considering only simplified problems.
As we show in the next sections, this obvious
observation has far reaching consequences, and can be used
to give simpler proofs for some results one finds in the nonlinear
optimization literature.

\section{The canonical form}
\label{secRank}
In this section we present a simple canonical form for the
classical nonlinear programming problem \pRef{nlp},
which can be used when
the derivative of the equality constraints has full
rank but the inequality constraints are degenerated.
In this canonical form the variables are
$\wvec{y} = \wlr{y_1,\dots,y_m}^{\wtr{}}$,
$\wvec{z} = \wlr{z_1,\dots,z_r}^{\wtr{}}$
and
$\wvec{w} = \wlr{w_1,\dots,z_{n - m - r}}^{\wtr{}}$,
where
\[
r := \wfc{\wrm{rank}}{\wfc{Dhg}{0}} - m.
\]
The equality constraints are given by
\[
\wfc{\hat{h}}{\wvec{y},\wvec{z},\wvec{w}} = \wvec{y} = 0,
\]
and there are two groups of inequality constraints.
The first one is given by
\[
\wfc{\hat{g}}{\wvec{y},\wvec{z},\wvec{w}} = \wvec{z} = 0.
\]
The second group of inequalities is
is given by a potentially
complicated function $c$:
\[
\wfc{c}{\wvec{y},\wvec{z},\wvec{w}} \leq 0 \in \wrn{p},
\]
about which, in principle, we know that
\pbDef{condc}
\wfc{c}{0,0,0} = 0 \hspace{1cm} \wrm{and} \hspace{1cm}
\wfc{D_w c}{0,0,0} = 0,
\peDef{condc}
but have no control over $\wfc{D_yc}{0,0,0}$ or $\wfc{D_z}{0,0,0}$.
The domain of these functions is a product of open sets
$\wcal{Y} \times \wcal{Z} \times \wcal{W}$ containing
$0$. In summary, we have the canonical nonlinear programming problem
\pbDef{cnlp}
\begin{array}{llcl}
\wrm{minimize}      &  \wfc{f}{\wvec{y},\wvec{z},\wvec{w}} &      & \\
\wrm{subject \ to } &  \wvec{y} & = & 0 \in \wrn{m},    \\
                    &  \wvec{z} & \leq & 0 \in \wrn{r}, \\
                    &  \wfc{c}{\wvec{x},\wvec{y},\wvec{w}} & \leq & 0 \in \wrn{p}, \\
                    & \wvec{y} \in \wcal{Y} \subset \wrn{m}, & \wvec{z} \in \wcal{Z} \subset \wrn{r} & \wvec{w} \in \wcal{W} \subset \wrn{n - m - r}.
\end{array}
\peDef{cnlp}
In this problem, the derivative of the constraints has the
simple form
\pbDef{dcanon}
\wfc{Dhgc}{\wvec{x},\wvec{y},\wvec{z}} =
\left(
\begin{array}{ccc}
\wvec{I}_{m \times m}  & 0                     & 0 \\
0                      & \wvec{I}_{r \times r} & 0 \\
 \wfc{D_yc}{\wvec{y},\wvec{z},\wvec{w}}  &
 \wfc{D_zc}{\wvec{y},\wvec{z},\wvec{w}}  & \wfc{D_wc}{\wvec{y},\wvec{z},\wvec{w}}
\end{array}
\right),
\peDef{dcanon}
with $\wfc{D_wc}{0,0,0}= 0$, where $\wvec{I}_{m \times m}$
is the $m \times m$ identity matrix.
The next theorem shows that nonlinear
programming problems can be reduced to the canonical form
\pRef{cnlp} under mild assumptions, and the discussion
in the previous section shows that the form of the
first and second order conditions and the Lagrange multipliers
do not change in this reduction. In many relevant situations
we can then use Theorem \ref{thmRank}
below and say rigorously
\begin{quote}
``without loss of generality, we can assume that our
linear programming problem is of the form \pRef{cnlp}''
\end{quote}
The purpose of the present article is to
call the readers attention to this simple intuitive idea,
which is formalized by Theorem \ref{thmRank} and
its proof,  which we now present.

\pbTheoremBT{thmRank}{The Canonical Form}
Suppose that the functions $h$ and $g$ in problem \pRef{nlp}
are of class $C^s$, with $s \geq 1$,
in a neighborhood $\wcal{A}$ of $0 \in \wrn{n}$,
and $\wfc{h}{0} = 0$ and $\wfc{g}{0} = 0$.
If the $n \times \wlr{m + r}$ matrix
\pbDef{chess}
\wlr{ \wgradf{h_1}{0}, \dots, \wgradf{h_m}{0},
\wgradf{g_1}{0}, \dots, \wgradf{g_r}{0}}
\peDef{chess}
has rank $m + r$ then there exist open sets $\wcal{Y} \subset \wrn{m}$,
$\wcal{Z} \subset \wrn{r}$ and $\wcal{W} \in \wrn{n - m - r}$,
with $\wlr{0,0,0} \in \wcal{Y} \times \wcal{Z} \times \wcal{W}$,
and a diffeomorphim
\[
q: \wcal{Y} \times \wcal{Z} \times \wcal{W}
\rightarrow \wfc{q}{\wcal{Y} \times \wcal{Z} \times \wcal{W}} \subset \wcal{A}
\]
of class $C^{s}$ with $\wfc{q}{0,0,0} = 0$
such that, for $1 \leq  \ell \leq r$,
\[
\wfc{\hat{h}}{\wvec{y},\wvec{z},\wvec{w}} := \wfc{h}{\wfc{q}{\wvec{y},\wvec{z},\wvec{w}}} = \wvec{y}
\hspace{1cm} \wrm{and} \hspace{1cm}
\wfc{\hat{g}_\ell}{\wvec{y},\wvec{z},\wvec{w}} := \wfc{g_\ell}{\wfc{q}{\wvec{y},\wvec{z},\wvec{w}}}
= z_\ell
\]
and for $1 \leq \ell \leq p$ the functions
\[
\wfc{c_\ell}{\wvec{y},\wvec{z},\wvec{w}} :=
\wfc{g_{r + \ell}}{\wfc{q}{\wvec{y},\wvec{z},\wvec{w}}}
\]
are such that $\wfc{c}{0,0,0} = 0$ and $\wfc{D_wc}{0,0,0}= 0$.
\peTheorem{thmRank}

Theorem \ref{thmRank} is a direct consequence of the 
Projection Lemma below, which is a corollary of
the Inverse Function Theorem (see Thm 2-13 in page 43 of \cite{Spivak}.)
Besides Lemma \ref{lemSubmersion}, we only need to note that
the equality $\wfc{D_w}{0,0,0} = 0$
follows from Equation \pRef{dcanon} and the
assumption that $\wfc{Dhg}{0}$ has rank $m + r$.

\pbLemmaBT{lemSubmersion}{The Projection Lemma}
Let $\wcal{A}$ be a neighborhood of $0 \in \wrn{n + k}$
and let $f: \wcal{A} \rightarrow \wrn{n}$ a function of class
$C^s$. If $\wfc{Df}{0}$
has rank $n$ then there exists a neighborhood
$\wcal{Y}$ of $0 \in \wrn{n}$,
a neighborhood $\wcal{Z}$ of $0 \in \wrn{k}$
and a diffeomorphism
\[
q: \wcal{Y} \times \wcal{Z} \rightarrow \wfc{q}{\wcal{Y} \times \wcal{Z}} \subset \wcal{A}
\]
of class $C^s$
such that $\wfc{f}{\wfc{q}{\wvec{y},\wvec{z}}} = \wvec{y}$.
\peLemma{lemSubmersion}

\section{A proof of Andreani's Conjecture}
\label{secNino}
In this section we use the canonical form
\pRef{cnlp} to prove the conjecture by
Andreani, Martinez and Schuverdt mentioned in the introduction.
This is an interesting application of the canonical form
because there were several failed attempts to find
a proof of this conjecture by other means.

For instance, the authors of \cite{Gabriel} attempted to obtain an appropriate
coordinate system, by using a version of the Singular Value
Decomposition, and succeeded in proving new particular cases
of Andreani's conjecture with this approach. However, they
did not prove the conjecture because their
decomposition is not as effective as the canonical form:
it has ``high order terms'' in places in which the
canonical decomposition
has exact zeros, and the technicalities required to handle
these terms precluded them from obtaining a
proof for which they had found all the other ingredients.
This shows that good choices of variables go beyond
controlling ``high order terms'': we actually want
to eliminate them, and Andreani's
conjecture is one of the fortunate cases in
which this is possible.

We now prove Andreani's Conjecture
using the arguments presented in the previous sections. Along the
proof we resort to two linear algebraic lemmas, which are proved
in the appendix. These lemmas are variations of results already
presented in other references \cite{Gabriel,GabrielB}, and we provide
their proofs to make the article self contained.
The first step to prove Theorem \ref{thmNino} is to use
Theorem \ref{thmRank} to
reduce problem \pRef{nlp} to the canonical form \pRef{cnlp}.
In order to that we note that the Mangasarian-Fromovitz
constraint qualification requires that
$\wfc{Dh}{0}$ has rank $m$ and recall that
we use $m + r$ to denote the rank of
the matrix $\wfc{Dhg}{0}$ in Equation \pRef{jhg}.
Therefore,
by changing the order of the inequality constraints if
necessary, we can assume that the matrix in Equation
\pRef{chess} in the statement of Theorem \ref{thmRank}
has rank $m + r$. We can then use this
theorem to analyze Andreani's conjecture 
and assume without loss of generality
that our linear programming problem has the form
\pRef{cnlp}, because in Section \ref{secChange} we have
shown that Equation \pRef{eqWeak} and the Mangasarian-Fromovitz
constraint qualification are invariant under
changes of coordinates.

Due to Equations \pRef{ortho} and \pRef{dcanon}, we can assume that the vector
$\wvec{d}$ in Equation \pRef{eqWeak} has the form
\[
\wvec{d}
=
\left(
\begin{array}{c}
0_{m} \\
0_{r} \\
\tilde{\wvec{d}}
\end{array}
\right)
\hspace{1cm} \wrm{with} \hspace{1cm}
\tilde{\wvec{d}} \in \wrn{p},
\]
and defining $\tilde{f}: \wcal{W} \rightarrow \wrone{}$
and  $\tilde{c}: \wcal{W} \rightarrow \wrone{}$ by
\[
\wfc{\tilde{f}}{\wvec{w}} := \wfc{f}{0_m,0_r,\wvec{w}}
\hspace{1cm} \wrm{and} \hspace{1cm}
\wfc{\tilde{c}}{\wvec{w}} := \wfc{c}{0_m,0_r,\wvec{w}},
\]
we can rewrite Equation \pRef{eqWeak} as
\pbDef{eqWeakT}
\tilde{\wvec{d}}^{\wtr{}}
\wlr{
\whessf{\tilde{f}}{0}  \ + \ \sum_{\ell = 1}^{p} \mu_{r + \ell} \whessf{\tilde{c}_\ell}{0}
}\tilde{\wvec{d}} \geq 0.
\peDef{eqWeakT}
Equation \pRef{dcanon} yields
\[
\wfc{\wrm{rank}}{\wfc{Dhgc}{\wvec{y},\wvec{z},\wvec{w}}} \leq m + r + 1
\
\Leftrightarrow
\
\wfc{\wrm{rank}}{\wfc{D_wc}{\wvec{y},\wvec{z},\wvec{w}}} \leq 1,
\]
and combining the observation that
\[
\wfc{D\tilde{c}}{\wvec{w}} = \wfc{D_wc}{0,0,\wvec{w}}
\]
with the next lemma we conclude that
there exist $\alpha_1,\dots,\alpha_p \in \wrone{}$ and
a symmetric matrix $\wvec{H} \in \wrn{p \times p}$ such that
\[
\whessf{\tilde{c}_\ell}{0} = \alpha_\ell \wvec{H}
\ \ \wrm{for} \ \ \ell = 1,\dots p.
\]

\pbLemmaBT{lemHessOne}{Hessians with rank at most one}
Let $\wcal{A}$ be a neighborhood of $0 \in \wrn{n}$,
and let $c_1,\dots,c_m$ be functions from
$\wcal{A}$ to $\wrone{}$ of class $C^2$. If
$\wfc{Dc}{0} = 0$ and $\wfc{\wrm{rank}}{\wfc{Dc}{\wvec{x}}} \leq 1$
for all $\wvec{x} \in \wcal{A}$ then there exist
$\alpha_1,\dots,\alpha_m \in \wrone{}$ and a
symmetric matrix
$\wvec{H} \in \wrn{n \times n}$ such that $\whessf{c_\ell}{0} = \alpha_\ell \wvec{H}$
for $\ell = 1,\dots,m$.
\peLemma{lemHessOne}

It follows that Equation \pRef{eqWeakT} is equivalent to
\pbDef{eqWeakTB}
\wfc{\tilde{S}}{\tilde{\wvec{d}},\gamma} := \tilde{\wvec{d}}^{\wtr{}}
\wlr{
\whessf{\tilde{f}}{0}  \ + \gamma \wvec{H}}
\tilde{\wvec{d}} \geq 0
\hspace{0.5cm} \wrm{for} \hspace{0.5cm}
\gamma = \sum_{\ell = 1}^{p} \alpha_\ell \mu_{r + \ell}.
\peDef{eqWeakTB}
The Mangasarian-Fromovitz
constraint qualification implies that for every
$\wvec{d}$ there exists $\lambda$ and $\mu > 0$
such that  $\wfc{S}{\wvec{d}} \geq 0$ in Equation \pRef{eqWeak}
(see \cite{Gabriel}.)
Since this equation is equivalent to Equation \pRef{eqWeakTB},
for every $\tilde{\wvec{d}}$ there exists $\gamma$
such that $\wfc{\tilde{S}}{\tilde{\wvec{d}}} \geq 0$ in Equation \pRef{eqWeakTB},
and in order to complete the proof we use the following lemma:

\pbLemmaBT{lemSemiSep}{Semidefinite separation}
Let $\wcal{I} \subset \wrone{}$ be a compact interval.
If the symmetric matrices $\wvec{A}, \wvec{B} \in \wrn{n \times n}$ are such that
for all $\wvec{x} \in \wrn{n}$ there exists
$\gamma_x \in \wcal{I}$ such that $\wvec{x}^{\wtr{}} \wlr{\wvec{A} + \gamma_x \wvec{B}} \wvec{x} \geq 0$
then there exists $\gamma^* \in \wcal{I}$ such that
\[
\wvec{x}^{\wtr{}} \wlr{\wvec{A} + \gamma^* \wvec{B}} \wvec{x} \geq 0
\]
for all $\wvec{x} \in \wrn{n}$.
\peLemma{lemSemiSep}
Lemma \ref{lemSemiSep} implies that there exists
\[
\gamma^*  = \sum_{\ell = 1}^p \alpha_\ell \mu^*_{r + \ell}
\]
such that $\wfc{\tilde{S}}{\tilde{\wvec{d}},\gamma^*} \geq 0$  in Equation \pRef{eqWeakTB}
for all $\tilde{\wvec{d}}$. The full set $\lambda^*$ and $\mu^*$ of Lagrange
multipliers containing these
$\mu^*_{r + 1}, \dots, \mu^*_{r + p}$
is as required by Andreani's Conjecture and we are done. \wqed{}

\appendix

\section{Linear Algebra}
\label{secLinAlg}
In this appendix we prove the results involving
Linear Algebra used in Section \ref{secNino}.
The proof of Lemma \ref{lemHessOne} is based on the
next two lemmas:

\pbLemmaBT{lemRankOneD}{Rank one columns}
If the symmetric matrices $\wvec{H}_1,\dots, \wvec{H}_m \in \wrn{n \times n}$
are such that the $n \times m$
matrix
\[
\wvec{A}_v := \wlr{\wvec{H}_j \wvec{v}, \dots, \wvec{H}_m \wvec{v}}
\]
has rank at most one for all $\wvec{v} \in \wrn{n}$ then there exist
$\alpha_1,\dots \alpha_m \in \wrone$ and a symmetric matrix
$\wvec{H} \in \wrn{n \times n} \setminus \wset{0}$
such that $\wvec{H}_j = \alpha_j \wvec{H}$ for $j = 1 \dots m$.
\peLemma{lemRankOneD}

\pbLemmaBT{lemHessLD}{Directional derivatives of rank one}
Let $\wcal{A}$ be a neighborhood of $0 \in \wrn{n}$
and let $h: \wcal{A} \rightarrow \wrn{n \times m}$ be a function
of class $C^1$, and for $1 \leq \ell \leq m$ let
$\wfc{h_\ell}{\wvec{x}}$ be the $\ell$th column of
$\wfc{h}{\wvec{x}}$.  If $\wfc{h}{0} = 0$
and $\wfc{\wrm{rank}}{\wfc{h}{\wvec{x}}} \leq 1$
for all $\wvec{x} \in \wcal{A}$ then
for every $\wvec{v} \in \wrn{n}$ the
$n \times m$ matrix
\[
\wvec{A}_v := \wlr{\wfc{Dh_1}{\wvec{x}} \wvec{v}, \dots,
                  \wfc{Dh_m}{\wvec{x}} \wvec{v}}
\]
has rank at most one.
\peLemma{lemHessLD}

The proof of Lemma \ref{lemSemiSep} uses the
following auxiliary results:

\pbTheoremBT{thmDines}{Dines' Theorem}
If the symmetric matrices $\wvec{A}, \wvec{B} \in \wrn{n \times n}$ are such that
for all $\wvec{x} \in \wrn{n} \setminus \wset{0}$
either $\wvec{x}^{\wtr{}} \wvec{A} \wvec{x} \neq 0$ or
$\wvec{x}^{\wtr{}} \wvec{B} \wvec{x} \neq 0$ then the set
\pbDef{rab}
\wfc{\wcal{R}}{\wvec{A},\wvec{B}} :=
\wset{ \wlr{\wvec{x}^{\wtr{}} \wvec{A} \wvec{x}, \, \wvec{x}^{\wtr{}} \wvec{B} \wvec{x}}, \ \
\wvec{x} \in \wrn{n}} \subset \wrn{2}
\peDef{rab}
is either $\wrn{2}$ itself or a closed angular sector of angle less than $\pi$.
\peTheorem{thmDines}

\pbLemmaBT{lemSep}{Definite separation}
Let $\wcal{I} \subset \wrone{}$ be a compact interval.
If the symmetric matrices $\wvec{A}, \wvec{B} \in \wrn{n \times n}$ are such that
for all $\wvec{x} \in \wrn{n} \setminus \wset{0}$ there exists
$\gamma_x \in \wcal{I}$ such that $\wvec{x}^{\wtr{}} \wlr{\wvec{A} + \gamma_x \wvec{B}} \wvec{x} > 0$
then there exists $\gamma^* \in \wcal{I}$ such that
\[
\wvec{x}^{\wtr{}} \wlr{\wvec{A} + \gamma^* \wvec{B}} \wvec{x} > 0
\]
for all $\wvec{x} \in \wrn{n} \setminus \wset{0}$.
\peLemma{lemSep}

Dines' Theorem is proved in \cite{Dines}, and in the rest of this
appendix we prove our lemmas in the order in which they were stated.

%%%%%%%%%%%%%%%%%%%%%%%%%%%%%%%%%%%%%%%%%%%%%%%%%%%%%%%%%%%%%%%%%%%%%%%%%%%%%%%%%%%%%
%
%%%%%%%%%%%%%%%%%%%%%%%%%%%%%%%%%%%%%%%%%%%%%%%%%%%%%%%%%%%%%%%%%%%%%%%%%%%%%%%%%%%%%

\pbProofB{Lemma}{lemHessOne}
Applying Lemma \ref{lemHessLD} to the
function $\wfc{h}{\wvec{x}} := \wfc{Dc}{\wvec{x}}$ we conclude that
for every $\wvec{v} \in \wrn{n}$ the $n \times m$ matrix
\[
\wvec{A}_v := \wlr{\whessf{c_1}{0} \wvec{v}, \dots,  \whessf{c_m}{0} \wvec{v}}
\]
has rank at most one, and Lemma \ref{lemRankOneD} yields the
coefficients $\alpha_j$ and the matrix $\wvec{H}$.
\peProof{Lemma}{lemHessOne}

%%%%%%%%%%%%%%%%%%%%%%%%%%%%%%%%%%%%%%%%%%%%%%%%%%%%%%%%%%%%%%%%%%%%%%%%%%%%%%%%%%%%%
%
%%%%%%%%%%%%%%%%%%%%%%%%%%%%%%%%%%%%%%%%%%%%%%%%%%%%%%%%%%%%%%%%%%%%%%%%%%%%%%%%%%%%%

\pbProofB{Lemma}{lemSemiSep}
For every $k \in \wn{}$, Lemma \ref{lemSep} yields $\gamma_k \in \wcal{I}$ such that
\[
\wvec{x}^{\wtr{}} \wlr{\wvec{A} + \gamma_k \wvec{B} + \frac{1}{k} \wvec{I}_{n \times n}} \wvec{x} \geq 0
\]
for all $\wvec{x} \in \wrn{n}$. Since the sequence $\gamma_k$ is bounded, it has a
subsequence which converges to some $\gamma^* \in \wcal{I}$. This $\gamma^*$ is as required
by Corollary \ref{lemSemiSep}.
\peProof{Lemma}{lemSemiSep}

%%%%%%%%%%%%%%%%%%%%%%%%%%%%%%%%%%%%%%%%%%%%%%%%%%%%%%%%%%%%%%%%%%%%%%%%%%%%%%%%%%%%%
%
%%%%%%%%%%%%%%%%%%%%%%%%%%%%%%%%%%%%%%%%%%%%%%%%%%%%%%%%%%%%%%%%%%%%%%%%%%%%%%%%%%%%%

\pbProofB{Lemma}{lemRankOneD}
Lemma \ref{lemRankOneD} holds when $n = 1$ or $m = 1$. Let us then assume that
it holds when for $n - 1 \geq 1$ or $m - 1 \geq 1$  and show that it also holds for
$m$ and $n$. If some $\wvec{H}_j$ is zero then we can
take $\alpha_j = 0$ and use induction for
\[
\wvec{H}_1,\dots,\wvec{H}_{j-1}, \wvec{H}_{j+1},\dots \wvec{H}_m.
\]
Therefore, we can assume that $\wvec{H}_j \neq 0$ for
all $j$. It follows that $\wvec{H}_1$ has an
eigenvalue decomposition
$\wvec{H}_1 = \wvec{Q} \wvec{D} \wvec{Q}^{\wtr{}}$ with
$d_{11} \neq 0$. By replacing $\wvec{H}_j$ by $\wvec{Q}^{\wtr{}} \wvec{H}_j \wvec{Q}$
for all $j$, we can assume that $\wvec{Q} = \wvec{I}_{n \times n}$.
Taking $\wvec{v} = 1/d_{11} \wvec{e}_1$
we obtain that $\wvec{H}_1 \wvec{v} =  \wvec{e}_1$ and the
hypothesis that the matrix $\wvec{A}_v$ has rank one implies that
$\wvec{H}_j \wvec{v} = \alpha_j \wvec{e}_1$ for all $j$.
Since the matrices $\wvec{H}_j$ are symmetric, all of them have the
form
\[
\wvec{H}_j =
\left(
\begin{array}{cc}
\alpha_j        & 0_{n -1}^{\wtr{}} \\
0_{n -1} & \wvec{H}'_j
\end{array}
\right),
\]
for symmetric matrices
$\wvec{H}_j' \in \wrn{\wlr{n-1} \times \wlr{n-1}} \setminus  \wset{0}$
which satisfy the hypothesis of Lemma \ref{lemRankOneD} with $n = n - 1$.
Since $\wvec{H}_1 \wvec{e}_1 = d_{11} \wvec{e}_1 \neq 0$, we have that
$\alpha_1 \neq 0$. Let then $\wvec{v}' \in \wrn{n -1}$ be such that
$\wvec{H}' \wvec{v}' \neq 0$, and write
\[
\wvec{v} =
\left(
\begin{array}{c}
1 \\
\wvec{v}'
\end{array}
\right)
\]
We have that
\[
\wvec{H}_1 \wvec{v} =
\alpha_1 \left(
\begin{array}{c}
 1 \\
\frac{\alpha_1'}{\alpha_1} \wvec{H}' \wvec{v}'
\end{array}
\right)
\hspace{0.5cm} \wrm{and} \hspace{0.5cm}
\wvec{H}_j \wvec{v} =
\left(
\begin{array}{c}
\alpha_j \\
\alpha_j' \wvec{H}' \wvec{v}'
\end{array}
\right)
\ \ \wrm{for} \ \ j > 1,
\]
and the $\wvec{H}_1 \wvec{v}$ and
$\wvec{H}_j \wvec{v}$ are aligned by hypothesis.
Moreover, either $\alpha_j \neq 0$ or $\alpha_j' \neq 0$,
because $\wvec{H}_j \neq 0$. It follows
that $\alpha_j' = \alpha_1' \alpha_j / \alpha_1$
for all $j$. As a result, $\wvec{H}_j = \alpha_j \wvec{H}$
where
\[
\wvec{H} :=
\left(
\begin{array}{cc}
1                       & 0^{\wtr{}}_{n-1} \\
0_{n-1} & \frac{\alpha_1'}{\alpha_1} \wvec{H}'
\end{array}
\right)
\]
and we are done.
\peProof{Lemma}{lemRankOneD}

%%%%%%%%%%%%%%%%%%%%%%%%%%%%%%%%%%%%%%%%%%%%%%%%%%%%%%%%%%%%%%%%%%%%%%%%%%%%%%%%%%%%%
%
%%%%%%%%%%%%%%%%%%%%%%%%%%%%%%%%%%%%%%%%%%%%%%%%%%%%%%%%%%%%%%%%%%%%%%%%%%%%%%%%%%%%%

\pbProofB{Lemma}{lemHessLD}
For every $\wvec{v} \in \wrn{n}$ and $1 \leq \ell \leq m$, the
facts that $h \in C^1$ and $\wfc{h}{0} = 0$
imply that
\[
\lim_{\delta \rightarrow 0}
\frac{1}{\delta} \wfc{h_\ell}{\delta \wvec{v}} = \wfc{Dh_\ell}{0} v  = \wvec{A}_v.
\]
Let $\rho > 0$ be
such that if $\wnorm{\wvec{B} - \wvec{A}_v} \leq \rho$ then
$\wfc{\wrm{rank}}{\wvec{B}} \geq \wfc{\wrm{rank}}{\wvec{A_v}}$.
Taking $\delta$ such that
\[
\wnorm{ \frac{1}{\delta} \wfc{h}{\delta \wvec{v}} - \wvec{A}_v } < \rho
\]
we obtain that
\[
1 \geq \wfc{\wrm{rank}}{\wfc{h}{\delta \wvec{v}}}
= \wfc{\wrm{rank}}{\frac{1}{\delta} \wfc{h}{\delta \wvec{v}}}
\geq \wfc{\wrm{rank}}{\wvec{A}_v}.
\]
Therefore, $\wfc{\wrm{rank}}{\wvec{A}_v} \leq 1$, and we are done.
\peProof{Lemma}{lemHessLD}

\pbProofB{Lemma}{lemSep}
Let us write $\wcal{I} = [a,b]$. If $a = b$ then
we could simply take $\gamma^* = a$ and we would be done.
Therefore, we can assume that $a < b$.
If $\wlr{u,v} \in  \wfc{\wcal{R}}{\wvec{A},\wvec{B}}$
then $ u + \gamma v \geq 0$ for some $\gamma \in [a,b]$. This implies that
\[
u \geq -\wabs{\gamma} \wabs{v}  \geq - \wlr{1 + \wabs{a} + \wabs{b}} \wabs{v}
\]
Therefore, $\wlr{- 2 \wlr{1 + \wabs{a} + \wabs{b}}, 1} \not \in \wfc{\wcal{R}}{\wvec{A},\wvec{B}}$
and  by Dines' Theorem $\wfc{\wcal{R}}{\wvec{A},\wvec{B}}$ is a closed pointed cone.
For each $\gamma \in \wcal{I}$ the set
\pbDef{ct}
\wcal{C}_\gamma := \wset{ \wlr{u, v} \in \wrn{2}, \ \ \wrm{with} \ u + \gamma v \leq 0}
\peDef{ct}
is a closed convex cone, and
\[
\wcal{C} := \bigcap_{\gamma \in \wcal{I}} \wcal{C}_\gamma = \bigcap_{\gamma \in [a,b]} \wcal{C}_\gamma
\]
is also a closed cone. Moreover, $\wcal{C}$ is pointed because $a < b$.
The hypothesis tells us that for every $\wlr{u,v} \in \wfc{\wcal{R}}{\wvec{A},\wvec{B}} \setminus
\wset{(0,0)}$ there exists $\gamma \in \wcal{I}$ such that $\wlr{u,v} \not \in \wcal{C}_\gamma$,
and this implies that $\wcal{C} \cap \wfc{\wcal{R}}{\wvec{A},\wvec{B}} = \wset{\wlr{0,0}}$.
Therefore, there exists $(c,d) \in \wrn{2}$ such that
\pbDef{sepcd}
c u + d v > 0 \hspace{0.1cm} \wrm{for} \hspace{0.1cm} (u,v) \in \wfc{\wcal{R}}{\wvec{A},\wvec{B}} \setminus \wset{\wlr{0,0}}
\peDef{sepcd}
and
\pbDef{sepcd2}
c u + d v < 0 \hspace{0.1cm} \wrm{for} \hspace{0.1cm} (u,v) \in \wcal{C} \setminus \wset{\wlr{0,0}}.
\peDef{sepcd2}
Equation \pRef{ct} shows that $(-1,0) \in \wcal{C}_\gamma$ for all $\gamma$. Therefore,
$(-1,0) \in \wcal{C}$ and Equation \pRef{sepcd2} implies that $c > 0$, and
by dividing Equations \pRef{sepcd} and \pRef{sepcd2} by $c$ we can assume that $c = 1$.
The point $\wlr{a,-1}$ belongs to $\wcal{C}_\gamma$ for all $\gamma \geq a$. Therefore,
$(a,-1) \in \wcal{C}$ and Equation \pRef{sepcd2} implies
that $a - d < 0$, that is $d > a$. Similarly,
The point $(-b,1)$ belongs to $\wcal{C}_\gamma$ for all $\gamma \leq b$, and
$-b + d < 0$, that is, $d < b$. In summary, we have that $d \in [a,b]$.
Finally, we take $\gamma^* = d$, and Equation \pRef{sepcd} with $c = 1$
shows that this is a valid choice.
\peProof{Lemma}{lemSep}

\bibliographystyle{amsplain}
\bibliography{nino}
\end{document}